\documentclass[graybox]{svmult}

\usepackage{type1cm}        % activate if the above 3 fonts are
\usepackage{makeidx}         % allows index generation
\usepackage{graphicx}        % standard LaTeX graphics tool
\usepackage{multicol}        % used for the two-column index
\usepackage[bottom]{footmisc}% places footnotes at page bottom

\usepackage{newtxtext}       % 
\usepackage[varvw]{newtxmath}       % selects Times Roman as basic font

\usepackage[utf8]{inputenc}

\usepackage{natbib}

\usepackage{mathrsfs}
\usepackage{graphicx}
\usepackage{subfigure}
\usepackage{color}

\DeclareMathOperator{\var}{var}

\DeclareMathOperator{\pr}{P}

\newcommand{\R}{\mathbb{R}}

\newcommand{\Ft}{\mathbf{F}_{\pm}}
 
\newcommand{\Fs}{\mathbf{F}_{\pm}^{(n)}}

\newcommand{\tZ}{\mathbf{Z}}
\newcommand{\tV}{\mathbf{V}}
\newcommand{\tW}{\mathbf{W}}

\newcommand{\tu}{\mathbf{u}}

\newcommand{\tT}{\mathbf{T}}
\newcommand{\tS}{\mathbf{S}}
\newcommand{\tX}{\mathbf{X}}

\newcommand{\tU}{\boldsymbol{U}}
\newcommand{\tz}{\boldsymbol{z}}
 
\newcommand{\Fsa}{\mathbf{F}_{\pm,a}^{(n)}}
\newcommand{\Rs}{R_{\pm}}
\newcommand{\Rsa}{R_{\pm,a}}
\newcommand{\Ssa}{\mathbf{S}_{\pm,a}}

\newcommand{\Zsa}{S_{\pm,a}}
\newcommand{\Zs}{S_{\pm}}

\renewcommand{\ts}{\boldsymbol{s}}

\newcommand{\In}{\mathrm{I}}

\renewcommand{\E}{\mathsf{E}}
\renewcommand{\phi}{\varphi}

\newcommand{\Dto}{\stackrel{D}{\to}}
\newcommand{\wto}{\stackrel{w}{\to}}
\newcommand{\Pto}{\stackrel{\pr}{\to}}

\newcommand{\calX}{\mathcal{X}}

\newcommand{\ta}{\mathbf{a}}
\newcommand{\eqD}{\stackrel{D}{=}}
\newcommand{\tJ}{\boldsymbol{J}}

\newcommand{\tmu}{\boldsymbol{\mu}}
\newcommand{\tg}{\mathbf{g}}
\newcommand{\tx}{\mathbf{x}}

\makeindex             % used for the subject index
\begin{document}

\title*{One-sample location tests based on center-outward signs and ranks}
\author{Daniel Hlubinka and \v{S}\'{a}rka Hudecov\'{a}}
\institute{Daniel Hlubinka \at Charles University, Faculty of Mathematics and Physics, Department of Probability and Statistics, \email{daniel.hlubinka@matfyz.cuni.cz}
\and Šárka Hudecová \at Charles University, Faculty of Mathematics and Physics, Department of Probability and Statistics, \email{hudecova@karlin.mff.cuni.cz}}
\maketitle

\abstract*{A multivariate one-sample location test based on the center-outward ranks and signs is considered, and two different testing procedures are proposed for centrally symmetric distributions. The first test is based on a random division of the data into two samples, while the second one uses a symmetrized sample. The asymptotic distributions of the proposed tests are provided. 
For univariate data, two variants of the symmetrized test statistic are shown to be equivalent to the standard sign and Wilcoxon test respectively. 
The small sample behavior of the proposed techniques is illustrated by a simulation study that also provides a power comparison for various transportation grids.}

\abstract{
A multivariate one-sample location test based on the center-outward ranks and signs is considered, and two different testing procedures are proposed for centrally symmetric distributions. The first test is based on a random division of the data into two samples, while the second one uses a symmetrized sample. The asymptotic distributions of the proposed tests are provided. 
For univariate data, two variants of the symmetrized test statistic are shown to be equivalent to the standard sign and Wilcoxon test respectively. 
The small sample behavior of the proposed techniques is illustrated by a simulation study that also provides a power comparison for various transportation grids.}

\keywords{center-outward ranks and signs, multivariate one-sample location test, distribution-freeness}

\section{Introduction}

Testing a location for a given sample is one of the basic problems in theoretical and applied statistics. It finds a number of applications, including a comparison of paired data. For univariate observations, the one-sample t-test is a basic technique which is valid asymptotically under a broad class of alternatives, but which may be far from optimal for distributions far from Gaussian. 
%A one-sample test of location is one of the most classical statistical
%techniques that are often used in the paired tests. 
There are several
well-known nonparametric tests of location, the sign test
and the one-sample Wilcoxon test being one of the most common, see \cite{hajek_book}. Both these tests use the notion of the ranks
and signs, which is, unfortunately, not directly extendable into a multivariate
setting. 

For a multivariate one-sample location problem, the classical procedure is the Hotelling's $T^2$ test, which is optimal
under Gaussian distribution but which can fail for non-elliptical distributions or distributions with heavy tails. 
There have been several attempts to define ranks for multivariate observations, see \cite{marc1} for an overview. \cite{OjaRan04}  proposed a one-sample location test based on spatial signs and spatial ranks, 
while \cite{randles} and \cite{OjaPain} studied tests defined via 
hyperplane based signs and ranks, see also recent computational improvements in \cite{HuSi}. The mentioned multivariate rank concepts require an additional assumption that the underlying distribution is elliptical, which can be limiting in applications.

Recently, \cite{cherno} and \cite{marc1} proposed a concept of multivariate ranks and signs 
derived from the theory of the optimal measure transportation. 
 The corresponding
 center-outward (c-o) ranks and signs have been applied in various multivariate statistical problems, see \citep{marc3,ghosal,shi, marcVAR} and further references therein. Among them \cite{marc2} proposed rank tests for
multiple-output regression, which involves a location two-sample problem as a special case. The main advantage of these c-o rank tests is that they are fully distribution free and asymptotically efficient.

The one-sample location test based on c-o ranks have not been studied in the literature. The regression framework from \cite{marc2} is not directly applicable without some  additional constraints, because the c-o ranks and signs are invariant with respect to a shift. 
Inspired by the univariate one-sample Wilcoxon test, one possible strategy is to assume some kind of symmetry of the underlying multivariate distribution and to derive a location test which employs this property. Two possible approaches are explored in this contribution for centrally symmetric multivariate distributions.

 %This contribution deals with a hypothesis about a location for a single centrally symmetric multivariate sample under. 
The paper is organized as follows. Section~\ref{sec:prel} defines the c-o ranks and signs.  Subsequently, two rank tests %based on this concept %on c-o ranks and signs 
 are proposed: The first one, discussed in Sections~\ref{sec:test1},  is based on a random division of the data into two samples and the two-sample location test from \cite{marc2}.
 The second approach, introduced in Section~\ref{sec:test2}, uses a symmetrized sample, and two test statistics are suggested. % and the asymptotic distribution of the test statistic is derived. 
The asymptotic distributions of the corresponding test statistics are derived.
Section~\ref{sec:d1} shows that for univariate data the tests are equivalent to the sign test and the Wilcoxon
test, respectively.  
 The finite sample behavior of the proposed methods is investigated in Monte Carlo study in Section~\ref{sec:simul}, which compares also the power for various transportation grids.

%We use the multivariate
%sign and ranks defined in \cite{marc1} to formulate and test two hypotheses on
%location. 
%Namely, we test the hypothesis on the centre of symmetry and
%a hypothesis on the MK median. Two of 

\section{Preliminaries}\label{sec:prel}

%The standard notation of %MT
%CO distribution function and the corresponding CO ranks and signs is used. %and quantile functions. 
Let $\tX$ be a $d$-dimensional random vector with
an absolutely continuous distribution $P_{\tX}$. % with a nonvanishing density.
 Let $\mathbb{B}^d$ be the closed unit ball in $\R^d$ and  $\mathcal{S}^d = \{\tz \in \mathbb{R}^d; \|\tz\| = 1\}$ be the unit sphere. 
  Denote $U_d$ the distribution of a random vector 
$\tU = R\tZ$ for $R$ and
  $\tZ$ independent, $R$ uniformly distributed on $[0,1]$, and $\tZ$
  uniformly distributed on
  $\mathcal{S}^d$. We say that a map $T$ pushes a distribution $P_1$ forward to $P_2$ if for $\tZ\sim P_1$ it holds that $T(\tZ)\sim P_2$. For the univariate case $d=1$, $U_1$ is simply the uniform distribution on the interval $[-1,1]$. 
 
 The center-outward distribution function (c-o df) $\Ft:\mathbb{R}^d \to \mathbb{B}^d$   of $\tX$ is defined as the gradient of a convex function pushing $P_{\tX}$ forward to $U_d$.
Under the existence of finite second moments $\Ft$ is the
$L_2$-optimal Monge-Kantorovich transportation, see \cite{marc1,marc3} for more details. 

\begin{remark}
The motivation for the definition of the c-o df $\Ft$ comes from the univariate situation. 
Let $X$ be a random variable with an absolutely continuous distribution with a cumulative distribution function (cdf) $F_X$. Then $X$ can be transformed monotonically to a random variable $W \sim U_1$, namely $W=2F_X(X) - 1=F_{\pm}(X) $. Moreover, if $\E X^2$ is finite, then the $L_2$ transport cost is minimal in a sense that
$\min_{V\sim U_1} \E (X-V)^2 = \E(X-F_{\pm}(X))^2$. 
%$\E(X-U)^2$ is minimal over all uniformly distributed random variables on $[-1,1]$. 
This property can be directly extended to an absolutely continuous random vector. The gradient of a convex function on $\mathbb{R}^d$ is a $d$-variate analogue of a  monotonic function on $\mathbb{R}$. It follows from a McCann Theorem \citep{mccann} that there is a $P_{\tX}$ almost surely unique gradient $\Ft$ of a convex function such that $\Ft(\tX)$ follows the uniform distribution $U_d$. Clearly, without the assumption of absolute continuity of $\tX$ there exists no transformation to uniformly distributed random variable even for $d=1$. %a univariate random variable$$. 
Assuming $\E\|\tX\|^2 < \infty$, the center-outward distribution function minimizes the mean square distance $\E\|\tX - \Ft(\tX)\|^2$ over all $U_d$ distributed random vectors.
\end{remark}

Let $\tX_1,\dots,\tX_n$  be a random sample from $P_{\tX}$. The empirical
  center-outward distribution function $\Fs$ is the mapping of
  $\tX_1,\dots,\tX_n$ to a regular grid $\mathcal{G}_n$ of $n$ points
  in the unit ball $\mathbb{B}^d$ such that 
  $\sum_{i=1}^{n}
  \|\Fs(\tX_i) - \tX_i\|$ is minimal. 
  In the following, the grid $\mathcal{G}_n$ is always assumed to satisfy the condition:  
  \[
  \text{The discrete uniform distribution on
  } \mathcal{G}_n \text{ converges weakly to } U_d \text{ as } n\to\infty.
 % \end{equation}
 \]
 One possibility is to take 
  \begin{equation}\label{eq:Gn}
  \mathcal{G}_n=\{\boldsymbol{g}_{ij}, i=1,\dots,n_R, j=1,\dots,n_S\}, \quad \text{ for } \boldsymbol{g}_{ij}= r_i \tu_j,
  \end{equation}
   where  $n=n_Rn_S$, $r_i=i/(n_R+1)$, and $\tu_j$ are  unit vectors, distributed approximately regularly over $\mathcal{S}^d$, or a union of $\mathcal{G}_n$ from \eqref{eq:Gn} and $n_0$ replications of $\{\boldsymbol{0}\}$; in that case $n=n_0+n_Rn_S$.   
For $d=2$ it is possible to construct a completely regular grid for any given $n_S$. For $d>2$  the unit vectors can  be obtained as a suitable transformation of a  low discrepancy sequence in $\R^{d-1}$, as Halton or Sobol, see \cite{fang} and our Section~\ref{sec:simul} for more details, where some alternative grids are considered as well.  If $n_0=1$, then $\bigl(\Fs\bigr)^{-1}(\boldsymbol{0})$ can be considered as an estimate of the sample location (sample median), see also discussion in \cite[Section 5]{marc1}.

%\textcolor{blue}{The center-outward distribution function also defines
%  a median (set). \citet{marc1} define the median as the intersection
%  $\bigcap_{0 < q < 1}\{\tx;\Ft(\tx)\leq q\}$. The empirical
%  counterpart is the point $\tx$ such that $\Fs(\tx)=\boldsymbol{0}$
%  if $\boldsymbol{0} \in \mathcal{G}_n$, or the intersection
%  $\bigcap_{0 < q < 1}\{\tx;\overline{\mathbf{F}}_{\pm}^{(n)}(\tx)\leq
%  q\}$, where $\overline{\mathbf{F}}_{\pm}^{(n)}$ is the smooth
%  interpolation of $\Fs$, see \citet{marc1}.}

 Finally, the center-outward rank $R_i$ and sign $\tS_i$ of
  $\tX_i$ are defined as
  \[
    R_i = (n_R+1)\|\Fs(\tX_i)\|, \quad \tS_i = \mathbb{I}_{[\Fs(\tX_i)\ne \boldsymbol{0}]}
    \frac{\Fs(\tX_i)}{\|\Fs(\tX_i)\|},
  \]
  so $\Fs(\tX_i) = R_i/(n_R+1)\cdot \tS_i$.

The center-outward distribution function
  $\Ft$ is invariant to a shift and it is equivariant to 
  orthogonal transformations. The empirical version $\Fs$ is
 also invariant to a shift and it is
   equivariant to an orthogonal transformation given by an orthogonal
  matrix $\mathbf{O}$ if the grid $\mathcal{G}_n$ is replaced by the
  transformed grid $\mathbf{O}\mathcal{G}_n$, see Proposition 2.2 in
  \citet{marc2}. In that case, the orthogonal transformation preserves the ranks and the angles among the signs. This is not the case for a general scaling transformation.

\section{Test of a center of symmetry}\label{sec:test}

Consider a $d$-dimensional random vector $\tX$ with a distribution which is centrally 
symmetric around a point $\tmu_S$, i.e.,
$\pr[\tX-\tmu_S \in A] = \pr[\tmu_S - \tX \in A]$ for any Borel set
$A \in \mathcal{B}(\R^d)$, see also \cite{serfling}. The hypothesis of interest is $H_0: \tmu_S = \tmu_0$ for some specified $\tmu_0\in\mathbb{R}^d$. It can be assumed without loss of generality that $\tmu_0=\boldsymbol{0}$, so 
%We need only consider the case $\mu_S=0$ as the
%general case can be always transformed to the special case by the
%shift $\tX-\mu_S$. Under this assumption we test the hypothesis
\[
 H_0: \tmu_S = \boldsymbol{0} \ \text{ against } \ H_1: \tmu_S \neq \boldsymbol{0}.
\]

%The point of symmetry $\tmu_S$ is preserved as the point of
%symmetry of the center-outward distribution function $\Ft$. 
%We use the fact in the construction of the tests of symmetry below.

The next proposition shows that under $H_0$ the center-outward distribution function is also centrally symmetric around zero. 

\begin{proposition}\label{prop1} Under $H_0$ it holds  that
%\[
$
\Ft(-\tx)=-\Ft(\tx),$ for all  $ \tx\in\R^d.
$
%\]
\end{proposition}
\begin{proof}
 % The equality follows from the symmetry of the distribution of $\tX$
 % and the construction of $\Ft$. % that $\Ft(\tx)=-\Ft(\tx)$. \textcolor{red}{Is it correct?}
  Recall  the  construction of $\Ft$ from \cite{marc1}: $\Ft$ is $\nabla \phi(\tx)$ for
  $\phi(\tx)=\sup_{\tu\in S^d}(\tu^\top\tx-\psi(\tu))$,
  where   $\psi(\tu)$ is the unique convex function such that $\nabla \psi$ 
 pushes the distribution $U_d$ forward to
  $P_{\tX}$ and $\psi(\boldsymbol{0})=0$. 
 It follows from the
  symmetry of both distributions that $-\nabla\psi(-\tu)$ also pushes
   $U_d$ forward to $P_{\tX}$. The uniqueness of $\psi$ gives
  $\psi(-\tu)=\psi(\tu)$. From this we conclude that
  $\phi(-\tx)=\sup_{\tu\in S^d}(-\tu^\top(-\tx)-\psi(-\tu))=\phi(\tx)$ and
  consequently
  $\Ft(\tx)=\nabla \phi(\tx) = -\nabla\phi(-\tx) =-\Ft(-\tx)$, which
  completes the proof.
\end{proof}

In the following consider a random sample $\mathcal{X}_n = (\tX_1,\dots,\tX_n)$ from a distribution which is centrally symmetric distribution around $\tmu_S$ and which has c-o df $\Ft$. 
Two  approaches for testing $H_0$   are proposed in the next two sections.
 %The first one, introduced in Section~\ref{sec:test1}, is 
%a variant of
%permutation test with observations multiplied by random
%signs.
%based on a two sample location test from \cite{marc2} applied to two randomly constructed sub-samples with random signs. 
% The second test, given in Section~\ref{sec:test2}, uses an augmented symmetrized sample. 

\subsection{Test with random signs}\label{sec:test1}

Assume that $H_0$ holds.  Then the distributions
of $\tX$ and $-\tX$ coincide, and, moreover, if
a random variable $Y$ is independent with $\tX$ and distributed uniformly on $\{-1,1\}$, that is
$\pr[Y=1] = \pr[Y=-1] = 1/2$, denoted as $\mathcal{U}\{-1,1\}$, then the random variables $\tX$ and
$\tX^{\circ} = Y\tX$ follow the same distribution. % and their center-outward distribution functions are the same.

Consider a %random sample $\mathcal{X}_n = (\tX_1,\dots,\tX_n)$ and a
sequence $Y_1,\dots,Y_n$ of iid $\mathcal{U}\{-1,1\}$ random variables independent with the
random sample $\mathcal{X}_n$.
Then the distribution of  $\mathcal{X}_n$ is the same as the distribution of
$\mathcal{X}_n^{\circ} = (\tX_{1}^{\circ}, \dots,\tX_{n}^{\circ})$, for $\tX_{i}^{\circ}=Y_i\tX_i$, 
and the center-outward distribution functions of $\mathcal{X}_n$ and
$\mathcal{X}_n^{\circ}$ are equal. In
particular, denoting $\mathcal{I}^{+} = \{i; Y_i = 1\}$,
$n_+ = |\mathcal{I}^{+}|$, and
\[
  \mathcal{X}^{+} = \{\tX_i; i \in \mathcal{I}^{+}\}, \quad
  \mathcal{X}^{-} = \{\tX_i; i \not\in \mathcal{I}^{+}\}
\]
it follows that $\mathcal{X}^{+}$ and $\mathcal{X}^{-}$ are two
independent samples from the same distribution. Note that the sample
size $n_{+}$ is random and follows the binomial distribution with
parameters $(n,1/2)$, which is the price we pay for the independence of
the two samples $\mathcal{X}^{+}$ and $\mathcal{X}^{-}$.

The main idea of the test with random signs is to test the hypothesis
$\mathrm{H}_0: \tX \stackrel{d}= -\tX$ using the two-sample location test  of \cite{marc2} on the independent samples $\mathcal{X}^{+}$
and $\mathcal{X}^{-}$. % Denote $\bigl(\tX^{(n)}_1,\dots,\tX^{(n)}_{n_1}\bigr) =\mathcal{X}^{+}$ and $\bigl(\tX^{(n)}_{n_1+1},\dots,\tX^{(n)}_{n}\bigr) =
The following proposition follows from \cite[Proposition 3.2]{marc2}.

%\mathcal{X}^{-}$. 

\begin{proposition}
Let $\Fs$ be the sample c-o distribution function computed for $\mathcal{X}_n^{\circ}$ and a grid $\mathcal{G}_n$, and let $R_i$ and $\tS_i$ be the corresponding rank and sign of $\tX_i^{\circ}$. 
Let $J:[0,1)\to\R$ be a score function continuous on $(0,1)$ with a bounded variation and $\int_0^1 J^2(u) \mathrm{d}u<\infty$, and define $\boldsymbol{J}(\tu) = \mathbb{I}_{\tu\ne \boldsymbol{0}} J(\|\tu\|)\cdot \tu/\|\tu\|$ for $\tu \in \mathbb{B}^d$. Denote as
\[
K_J (\mathcal{G}_n)=\sum_{\boldsymbol{g}\in\mathcal{G}_n} \boldsymbol{J}(\boldsymbol{g}). 
\]
Then   
\[
\boldsymbol{T}=\sqrt{\frac{n}{n_+(n-n_+)}}\left[\sum_{i \in \mathcal{I}^{+}} J\left(\frac{R_i}{n_R+1}\right) \tS_i - \frac{n_+}{n} K_J (\mathcal{G}_n) \right]
\]
converges in distribution for %$n_R,n_S\to\infty$ 
$n\to\infty$ to a centred normal distribution with variance matrix $(1/d)\cdot \int_{0}^1 J(u) \mathrm{d}u \cdot \boldsymbol{I}_d$ with probability one, where $\boldsymbol{I}_d$ is the identity matrix.
\end{proposition}

\begin{remark}
 The convergence in distribution holds with probability one since $n_{+}$ is a random variable. According to \cite{marc2}  $n_{+}/n$ must converge to a constant $c \in (0,1)$, as $n \to \infty$, to ensure asymptotic normality. Since the a.s. convergence of $n_{+}/n$ to $1/2$ follows from the strong law of large numbers and construction of $\mathcal{X}^{\circ}_{n}$ we get the asymptotic normality with probability one along the sequence of realizations of $\{Y_i\}$.
\end{remark}

%The hypothesis $H_0$ is rejected if
%\[
%Q = \frac{dn}{n_+(n-n_+) \int_0^1 J^2(u)du} \| \boldsymbol{T}\|^2
%\]

The hypothesis $H_0$ is rejected if
\[
Q = \frac{d}{ \int_0^1 J^2(u)\mathrm{d}u} \| \boldsymbol{T}\|^2
\]
exceeds $1-\alpha$ quantile of $\chi_d^2$ distribution.

Note that if the null hypothesis $H_0$ is violated, then the distributions of $\tX$ and $\tX^{\circ}$ differ, i.e., $\mathcal{X}^{+}$ and $\mathcal{X}^{-}$ are not
sampled from the same distribution. Under $H_1$, the distribution of $\tX$
is not symmetric with respect to the origin, while the distribution of
$\tX^{\circ}$ is  symmetric around $\boldsymbol{0}$ regardless the
validity of $H_0$. Monte Carlo simulations presented in Section~\ref{sec:simul} illustrate that in such case the test reveals the difference in the locations of the two samples  $\mathcal{X}^{+}$ and  $\mathcal{X}^{-}$.

\begin{remark}
The asymptotic distribution of the test statistic $Q$ under a sequence of local alternatives can be derived from the two-sample situation studied in \cite[Proposition 5.1 and Section~5.3.1]{marc2}. We present the result for the special case of spherically symmetric  distributions. 

Let the distribution of $\tX$ be spherically symmetric with a radial density $g$, i.e. the density $f$ of $\tX$ is proportional to $g(\sqrt{\tx^\top \tx})$. Denote as $f_r$ and $F_r$ the density and cumulative distribution function  of $\|\tX\|$ respectively.  
Assume that $g$ is mean square differentiable and $\varphi_g= - (g^{1/2})'/g^{1/2}$ is such that $\int_0^1 \varphi_g(F_r^{-1}(u))\mathrm{d} u <\infty$. 
Let $\boldsymbol{h}\in\R^d$ be some non-zero vector.
Under the sequence of local alternatives $H_1: \boldsymbol{\theta} = \boldsymbol{h}/\sqrt{n}$,  the test statistic $Q$ has asymptotically for $n\to\infty$ a non-central $\chi^2_d$ distribution with a non-centrality parameter
\[
q=\frac{\left[\int_0^1 J(u)\phi_g\left(F_r^{-1}(u)\right) \mathrm{d}u\right]^2}{d \int_0^1 J^2(u)\mathrm{d}u} \boldsymbol{h}^\top\boldsymbol{h}. 
\]
It easily follows from the latter formula that the test based on the score function $J_0 = \phi_g\circ F_r^{-1}$ is locally asymptotically optimal, c.f. \cite[Corollary 5.2 (ii)]{marc2}. Moreover, the asymptotic relative efficiency (ARE) of the test with respect to the Hotelling's $T^2$ test is given as
\[
\mathrm{ARE} =  \frac{\left[\int_0^1 J(u)\phi_g\left(F_r^{-1}(u)\right) \mathrm{d}u\right]^2  \cdot \E \|\tX\|^2}{d^2 \int_0^1 J^2(u)\mathrm{d}u}. 
\]
If $J=J_0$ then the ARE simplifies to $\int_0^1 J_0^2(u)\mathrm{d}u \cdot \E \|\tX\|^2/d^2$. 

If $\tX$ follows the normal distribution $\mathsf{N}(\boldsymbol{\theta},\boldsymbol{I}_d)$ then the optimal score function corresponds to the van der Waerden scores $J_0(u)=\sqrt{G_d^{-1}(u)}$, where $G_d$ is the cdf of the $\chi^2_d$ distribution, and in that case ARE equals $1$. If $J(u)=u$ (Wilcoxon scores) then straightforward calculations give that
\[
\mathrm{ARE}=\frac{3}{d}\left( \int_0^{\infty} \sqrt{y} G_d(y) \mathrm{d} G_d(y)\right)^2. 
\]
This for $d=1$ gives $\mathrm{ARE}=3/\pi$, which is the well-known ARE of the one sample Wilcoxon test with respect to the t-test. For $d=2$, we get $\mathrm{ARE} \doteq 0.985$ and the values decrease with increasing dimension $d$. For instance for $d=10$, $\mathrm{ARE}=0.907$. 
\end{remark}

\subsection{Symmetrized sample}
\label{sec:test2}
%Let $\tX_1,\dots,\tX_n$ still be a random sample from the distribution with CO distribution function $\Ft$. 

Consider now the set of size $2n$ of observations
$\calX_a=\{\tX_1,\dots,\tX_n,-\tX_1,\dots,-\tX_n\}$. Recall that under
the null hypothesis, $-\tX_1$ has the same distribution as $\tX_1$, so all the random vectors from $\calX_a$ have the same distribution under $H_0$.
Let $\Fsa$ be the sample center-outward distribution function  computed from $\calX_a$ which takes values in a  grid
  $\mathcal{G}_{2n}$ of $2n$ points satisfying
 \begin{equation}\label{eq:Gsym}
  \tg\in \mathcal{G}_{2n} \Rightarrow -\tg\in\mathcal{G}_{2n}. 
  \end{equation}
Hence, $\Fsa(\tX_i)$ stands for the value
corresponding to $\tX_i$, and  $\Rsa(\tX_i)$ and  $\Ssa(\tX_i)$ are its sample
c-o rank and sign respectively.  This slightly different notation is used in order to distinguish these quantities from those introduced in Section~\ref{sec:prel}. 
%Furthermore, when writing $n\to\infty$, we assume that both $n_R$ and $n_S$ converge to $\infty$.

\begin{proposition} %Let $H_0$ hold and let $\Fsa$ take values in the grid
%  $\mathcal{G}$ of $2n$ points, which satisfies
  %$\tg\in \mathcal{G} \Rightarrow -\tg\in\mathcal{G}$. Then
  Under $H_0$
  \[
    \Fsa(X_i)=- \Fsa(-X_i).
  \]
\end{proposition}
\begin{proof}
  Recall that $\Fsa$ is a mapping $g:\calX_a\to \mathcal{G}_{2n}$, which
  minimizes
\[
\sum_{i=1}^n \|\tX_i-g(\tX_i)\| + \sum_{i=1}^n \|-\tX_i-g(-\tX_i)\|^2, 
\]
or equivalently, maximizes
\begin{equation}\label{eq:equiv.cond.F}
\sum_{i=1}^n \tX_i^{\top}[g(\tX_i)-g(-\tX_i)].
\end{equation}
To prove the claim let us first consider a norm ordering of the sample
$\tX_1,\dots, \tX_n$ such that
\[
  \|\tX_{1:n}\| \geq  \|\tX_{2:n}\| \geq \dots \geq  \|\tX_{n:n}\|,
\]
and note that the inequalities are a.s. strict since the underlying
distribution of $\tX$ is absolutely continuous.

Recall that for two ordered sets of real numbers
$a_1 < a_2 < \dots < a_n$, and $b_1 < b_2 < \dots < b_n$ it holds that
$\sum a_i b_{\pi(i)} \leq \sum a_i b_i$ and the inequality is strict if
$\pi(i) \neq i$ for some $i$. Similar inequality holds for the sum of
inner products \eqref{eq:equiv.cond.F}.  The sum is maximized if we
match the ``large'' values of $\|\tX\|$ with $\tg \in \mathcal{G}_{2n}$ with as
most similar direction as possible. In other words, find $\tg \in \mathcal{G}_{2n}$ such
that the inner product $\langle \tX_{1:n}, \tg \rangle$ is maximal. Such
$\tg$ is unique a.s. (see the remark below). Set $g(\tX_{1:n}) = \tg =
-g(-\tX_{1:n})$. Clearly $\langle \tX_{1:n}, 2\tg \rangle > \langle \tX_{1:n}, \tg-\boldsymbol{h}
\rangle >$ for any $\boldsymbol{h} \neq -\tg$. Remove $\tX_{1:n}$ from the sample and
$\tg$ from the grid and continue with $\tX_{2:n}$ up to $\tX_{n:n}$. Hence,
with probability 1 we obtain $\Fsa(\tX_i)=- \Fsa(-\tX_i)$.

%\textcolor{blue}{Can this be used for the optimisation?}
%\textcolor{red}{To be added.}
\end{proof}

\begin{remark}
If there are two grid points $\tg$ and $\tg'$ such that
$\langle \tX_{\ell:n}, \tg \rangle = \langle \tX_{\ell:n}, \tg' \rangle$ then the
choice between $\tg$ and $\tg'$ depends on the sample points $\tX_{i:n}$, $i = \ell+1,\dots n$ with
probability 1. Hence, there is a sample point $\tX_{i:n}$ associated in the sense of the proof above with (without loss of generality) $\tg'$ such that $\langle \tX_{i:n}, \tg'
\rangle > \langle \tX_{i:n}, \tg \rangle$. Since
\[
  \langle \tX_{\ell:n}, 2\tg \rangle + \langle \tX_{i:n}, 2\tg' \rangle >
  \langle \tX_{\ell:n}, \tg + \tg' \rangle + \langle \tX_{i:n}, \tg'+\tg \rangle
\]
it becomes clear that the maximum of \eqref{eq:equiv.cond.F} is
attained if $\Fsa(\tX_i)=- \Fsa(-\tX_i)$. 
\end{remark}

\begin{proposition}\label{prop2}
%Let $\Fsa$ be as in Proposition~\ref{prop1}. Then 
Under $H_0$ 
$
\Fsa(\tX_1) \Pto \Ft(\tX_1)
$ as $n\to\infty$.
\end{proposition}

\begin{proof}
The proof goes along the lines of proof of Proposition 3.3 in \cite[Appendix F.3]{marc1} if we prove the weak convergence of the empirical distribution $\widehat{P}_{2n,a}$ of $\mathcal{X}_a$ to the distribution $P_{\tX}$. The latter convergence follows by Portmanteau theorem as for any closed set $F$ it holds
\[
\begin{split}
  \limsup_{n \to \infty}\widehat{P}_{2n,a}(F) & =  \limsup_{n\to\infty}\left(\frac{1}{2n}\sum_{i=1}^{n} \mathbb{I}_{[\tX_i \in F]} + \frac{1}{2n}\sum_{i=1}^{n} \mathbb{I}_{[-\tX_i \in F]}\right) \\
  & \leq \frac{1}{2}P_{\tX}(F) + \frac{1}{2}P_{-\tX}(F)
\end{split}
\]
 since the empirical distributions of $(\tX_1,\dots,\tX_n)$ and $(-\tX_1,\dots,-\tX_n)$ converge weakly to $P_{\tX}$ and $P_{-\tX}$, respectivelly.
 As $P_{\tX}(F) = P_{-\tX}(F)$ under the null hypothesis we get
for any closed set $F$ the inequality $\limsup_{n \to \infty}\widehat{P}_{2n,a}(F) \leq P_{\tX}(F)$ that implies $\widehat{P}_{2n,a} \wto P_{\tX}$.
\end{proof}

Recall that the grid $\mathcal{G}_{2n}$ is assumed to satisfy the condition \eqref{eq:Gsym}, which implies that
$\sum_{\tg\in\mathcal{G}_{2n}} \tg=\boldsymbol{0}$.  Consider two test statistics:
 \begin{align*}
\tT_{S,a}&=\frac{1}{\sqrt{n}}\sum_{i=1}^n \Ssa(\tX_i), \\
\tT_{F,a}&=\frac{1}{\sqrt{n}}\sum_{i=1}^n \Fsa(\tX_i) =\frac{1}{\sqrt{n}}\sum_{i=1}^n \Rsa(\tX_i)  \Ssa(\tX_i).
\end{align*}
The test statistic $\tT_{S,a}$ is based on signs only, while $\tT_{F,a}$ is a Wilcoxon type statistic.

\begin{proposition}\label{prop_as} Under $H_0$ %(and the assumption of  symmetry and the property of $\mathcal{G}$) 
\[
  d \| \tT_{S,a} \|^2 = \frac{d}{n} \left\| \sum_{i=1}^n
    \Ssa(\tX_i)\right\|^2\Dto \chi^2_d,
\]
and
\[
  3d \| \tT_{F,a} \| = \frac{3d}{n} \left\| \sum_{i=1}^n
    \Fsa(\tX_i)\right\|^2\Dto \chi^2_d
\]
as $n\to\infty$.
\end{proposition}

\begin{proof}
  Both statistics are special cases of
\begin{equation}\label{eq:Ta}
  \tT_{a}= \frac{1}{\sqrt{n}}\sum_{i=1}^n \boldsymbol{J}(\Fsa(\tX_i)) = \frac{1}{\sqrt{n}}\sum_{i=1}^n J(\|\Fsa(\tX_i)\|)\frac{\Fsa(\tX_i)}{\|\Fsa(\tX_i)\|},
\end{equation}
where the score function $J$ is chosen $J(r)=1$ for $\tT_{S,a}$ and $J(r)=r$ for
$\tT_{F,a}$. In both cases $\sum_{\tg\in\mathcal{G}_{2n}} \boldsymbol{J}(\tg)=\sum_{\tg\in\mathcal{G}_{2n}} J(\|\tg\|)\tg/\|\tg\| =\boldsymbol{0}$. In the following, we restrict $J$ to the two functions above.  
It follows from the property of the grid $\mathcal{G}_{2n}$ that
\[
  \frac{1}{\sqrt{2n}}\left(\sum_{i=1}^n \boldsymbol{J}(\Fsa(\tX_i)) - \sum_{i=1}^n
    \boldsymbol{J}(\Fsa(-\tX_i))\right) = \frac{2}{\sqrt{2n}} \sum_{i=1}^n
  \boldsymbol{J}(\Fsa(\tX_i)) - \boldsymbol{0} = \sqrt{2} \tT_{a}.
\]
Define
\[
  \tT^{e}_{a}= \frac{1}{\sqrt{2n}}\left( \sum_{i=1}^n \boldsymbol{J}(\Ft(\tX_i)) -
    \sum_{i=1}^n \boldsymbol{J}(\Ft(-\tX_i)) \right).
\]
Then, using Proposition~\ref{prop1} and since $\boldsymbol{J}(\tx)=-\boldsymbol{J}(-\tx)$,
\[
  \tT^{e}_{a}= \frac{\sqrt{2}}{\sqrt{n}} \sum_{i=1}^n \boldsymbol{J}(\Ft(\tX_i))
\]
and
\[
  \sqrt{2}\tT_{a}-\tT^{e}_{a} = \frac{1}{\sqrt{2n}}\sum_{i=1}^n (\ta_i
  -\ta^{-}_{i}),
\]
where $\ta_i = \boldsymbol{J}(\Fsa(\tX_i)) - \boldsymbol{J}(\Ft(\tX_i))$ and
$\ta^{-}_{i} = \boldsymbol{J}(\Fsa(-\tX_i)) - \boldsymbol{J}(\Ft(-\tX_i))=-\ta_i$ for
$i=1,\dots,n$. Then
\[
  \E \| \sqrt{2}\tT_{a}-\tT^{e}_{a} \|^2=\frac{1}{2n}4\sum_{i=1}^n\E
  \|\ta_i\|^2 + \frac{1}{2n}\sum\sum_{j\ne i} \E(\ta_i-\ta^{-}_{i})^{\top}(\ta_j-\ta^{-}_{j}).
\]
Because $\tX_1,\dots,\tX_n$ have the same distribution, then
$\E\|\ta_i\|^2=\E\|\ta_1\|^2$ for all $i=1,\dots,n$. Furthermore,
\[
  \E(\ta_i-\ta^{-}_{i})^{\top}(\ta_j-\ta^{-}_{j}) = \E\ta_i^{\top}\ta_j
  -\E\ta^{-\top}_{i} \ta_j-\E \ta_{i}^{\top}\ta^{-}_{j} +\E \ta^{-\top}_{i}\ta^{-}_{j}=0.
\]
The last equality follows from the fact that $(\tX_i,\tX_j) \eqD(-\tX_i,\tX_j)$, hence
$(\ta_i,\ta_{j}) \eqD (\ta^{-}_{i},\ta_j)$ and thus, all the expectations are the
same. Hence, it suffices to show that
\[
  \E \| \ta_1\|^2 = \E \| \boldsymbol{J}(\Fsa(\tX_i)) - \boldsymbol{J}(\Ft(\tX_i))\|^2 \to 0.
\]
Recall the convergence
$
  \Fsa(\tX_i)-\Ft(\tX_i) \Pto \boldsymbol{0}.
$
proved in Proposition~\ref{prop2}. Then it follows by the continuity of $J$ that 
$\tJ(\Fsa(\tX_1)) - \tJ(\Ft(\tX_1)) \Pto \boldsymbol{0}$ and
$\E \| \tJ(\Fsa(\tX_1)) \|^2\to \E \| \tJ(\Ft(\tX_1)) \|^{2}$. This shows that
$\sqrt{2}\tT_{a}-\tT^{e}_{a} =o_{q.m.}(1)$, and thus, the asymptotic
distribution of $\tT_{a}$ is the same as asymptotic distribution of
$\tT^{e}_{a}/\sqrt{2}$.  Recall that $\Ft(\tX_i) $ has the same
distribution as $\tV=U \cdot \tW$, where $U$ is uniformly distributed
on $[0,1]$ and $\tW$ is uniformly distributed on the unit sphere in
$\mathbb{R}^{d}$, and $U$ and $\tW$ are independent.  Then $\E \tV=0$ and
$\var \tV=\E U^2 \E\tW\tW'=\frac{1}{3} \var \tW =
\frac{1}{3d}\boldsymbol{I}_d$. The asymptotic distribution of $\tT^{e}_{a}$ follows from
the central limit theorem.

If $J(r)=r$ then 
\[
  \sqrt{2}\tT^{e}_{a}=\frac{1}{\sqrt{n}}\sum_{i=1}^n \Ft(\tX_i) \Dto
  \mathsf{N}_d\left(\boldsymbol{0},\frac{1}{3d}\boldsymbol{I}_d\right)
\]
and thus $ 3d \| \tT_{F,a} \| \Dto \chi^2_d$.
%\[
 % 3d \| \tT_{S,a} \| = \frac{3d}{n} \left\| \sum_{i=1}^n
  %  \Fsa(X_i)\right\|^2\Dto \chi^2_d.
%\]
If $J(r)=1$ then
 \[
   \sqrt{2}\tT^{e}_{a}= \frac{1}{\sqrt{n}}\sum_{i=1}^n
    \Ssa(\tX_i) \Dto \mathsf{N}_d\left(\boldsymbol{0},\frac{1}{d}\boldsymbol{I}_d\right)
\]
and thus $ d \| \tT_{S,a} \|^2  \Dto \chi^2_d$. 
%\[
 % d \| \tT_{S,a} \|^2 = \frac{d}{n} \left\| \sum_{i=1}^n
  %  \Ssa(X_i)\right\|^2\Dto \chi^2_d.
%\]
\end{proof}

\begin{remark}
The proof reveals that one could consider the general test statistic $\boldsymbol{T}_a$ defined in \eqref{eq:Ta} via a general score function $J$ satisfying assumptions of Proposition~2 such that $J(0)=0$. The discussion about the optimal choice of $J$ and the AREs of the test would be analogous to the one provided in Remark~3. We rather show in the next Section~\ref{sec:d1} that the tests based on the two test statistics from Proposition~5 are for $d=1$ asymptotically equivalent to the sign and one-sample Wilcon test respectively. 
\end{remark}

\subsection{One-dimensional analogue of the symmetrized test}\label{sec:d1}

Let $X_1,\dots,X_n$ be a random sample from a univariate continuous distribution, and 
consider the two test statistics proposed in Section~\ref{sec:test2}. % for univariate observations, i.e. for $d=1$. 

\begin{proposition}\label{prop6}
Let $R_1,\dots,R_n$ be the univariate ranks of
$|X_1|, \dots,|X_n|$. Then
\[
\Rsa(X_i)=R_i= \Rsa(-X_i), \quad  \Zsa(X_i)=(\In[X_i>0]-\In[X_i<0]) = -\Zsa(-X_i)
\]
and
\[
\Fsa(X_i)=\frac{R_i}{n+1}\left(\In[X_i>0]-\In[X_i<0]\right). 
\]
\end{proposition}
\begin{proof}
It is shown in Appendix B of \cite{marc1} that  if  $Z_1,\dots,Z_N$ is a random sample from a continuous univariate distribution and $R_1^Z,\dots,R_N^Z$ are the classical univariate ranks of $Z_1,\dots,Z_N$, then 
for even $N$
\begin{equation}\label{eq:ZR}
\Rs(Z_i)= \left|R_i^Z-\frac{N+1}{2}\right|+\frac12, \quad \Zs(Z_i)=\In\left[R_i^Z>\frac{N+1}{2}\right] - \In\left[R_i^Z<\frac{N+1}{2}\right],
\end{equation}
 and $\Fs(Z_i)=\Rs(Z_i) \Zs(Z_i) \left(N/2 +1\right)^{-1}$.

Let $R_{a}(X_i)$ stand for the (classical univariate) rank of $X_i$ in the augmented sample  $\calX_a = (X_1,\dots,X_n,-X_1,\dots,-X_n)$.
It is easy to see that 
$
R_{a}(X_i)+R_{a}(-X_i)=2n+1
$
with
\[
R_{a}(X_i)=\begin{cases} 
n+R_i & X_i>0,\\
n-R_i+1& X_i<0.
\end{cases}
\]
Moreover, $R_{a}(X_i)>n+1/2$ if and only if $X_i>0$.  Then it follows from \eqref{eq:ZR} that for $ X_i>0$ we have
\[
\Rsa(X_i)=\left| R_{a}(X_i)-\frac{2n+1}{2}\right|+\frac12 = \left| n+R_i-\frac{2n+1}{2}\right|+\frac12 = R_i = \Rsa(-X_i),
\]
and the same equality is shown analogously for $X_i<0$. 
Further, 
\[
\Zsa(X_i)=\In[X_i>0]-\In[X_i<0] = -\Zsa(-X_i),
\]
which finishes the proof.
\end{proof}

Since the distribution of $X_i$ is continuous, we can assume without loss of generality that $X_i\ne 0$ for all $i=1,\dots,n$. It follows from Proposition~6 that
%For the test statistic $T_{S,a}$, it follows that 
\[
\sqrt{n}T_{S,a}=\sum_{i=1}^n (\In[X_i>0]-\In[X_i<0])= 2\sum_{i=1}^n
\In[X_i>0]-n = 2W_S-n,
\]
 where $W_S = \sum_{i=1}^n \In[X_i>0]$ is the
test statistic of the univariate sign test. Hence, for $d=1$ the test based on
$T_{S,a}$ is exactly the univariate sign test. Note that it follows from the well known properties of $W_S$ that that $T_{S,a} \Dto \mathsf{N}(0,1)$, and
$T_{S,a}^2\Dto\chi_1^2$, which  is exactly what is claimed in
Proposition~\ref{prop_as} for $d=1$.
 
Similarly, for the Wilcoxon test statistic $T_{F,a}$, it follows that 
\begin{align*}
  \sqrt{n}T_{F,a}=&\frac{1}{n+1}\left( \sum_{i:X_i>0} R_i -
    \sum_{i:X_i<0} R_i \right) = \frac{1}{n+1}\left(2\sum_{i:X_i>0}
    R_i - \frac{n(n+1)}{2}\right)\\
     = &\frac{2W_W-\frac{n(n+1)}{2}}{n+1},
\end{align*}
where $W_W= \sum_{i:X_i>0} R_i $ is the test statistic of the one
sample Wilcoxon test. It is well known that
\[
  \frac{2W_W-\frac{n(n+1)}{2}}{\sqrt{1/3 n^3}} =
  \frac{\sqrt{3}}{\sqrt{n}}\frac{2W_W-\frac{n(n+1)}{2}}{n} \Dto
  \mathsf{N}(0,1).
\]
Hence, for $d=1$ the test based on $T_{F,a}$ is equivalent to the one-sample
Wilcoxon test, $\sqrt{3}T_{F,a} \Dto \mathsf{N}(0,1)$ and
$ 3\| T_{F,a}\|^2\Dto \chi^2_1.  $ This again coincides with the
asymptotics claimed in Proposition~\ref{prop_as} for $d=1$.

\section{Simulations}
\label{sec:simul}

 The performance of  the two proposed approaches % to one-sample center-outward test
  is explored in a Monte Carlo study conducted in R program,  \cite{R}. 
The power of each test is studied with respect to various distributional alternatives, and also with respect to various constructions of the grid for the empirical c-o distribution function. The following $d$-dimensional distributions were considered under the null hypothesis: 
\begin{enumerate}
\item[(a)] multivariate centered normal distribution with identity variance matrix,
\item[(b)] multivariate $t$ distribution with $df=1$ degrees of freedom, and identity scale matrix,  % X=Y/sqrt(U/df), where Y and U indep. and Y N(0,I) and U chi^2_df
%\item mixture of uniform distributions on $[0,1]^d$ and $[-1,0]^d$ with equal weights,  
\item[(c)] mixture of $\mathcal{L}(\tX)$ and $\mathcal{L}(-\tX)$ with equal weights, where $\tX$ is a $d$-dimensional vector with independent components with marginal exponential distribution with mean $1$.
\end{enumerate}
The multivariate Cauchy $t_1$ distribution in (b) is used for its heavy tails, while the distribution (c), in the following referred to as ``double-exponential'', represents a non-elliptical distribution. 
Notice that multivariate $t$-distributions with $df>1$ were also considered, and the corresponding results  are somehow in between the results obtained for (a) and (b), % $t_1$ and for the normal distribution in (a), 
so they are not presented in detail. 

Under the alternative, the above distributions were shifted by the vector $\delta \cdot \ts$, where $\delta\in\{0.0,\, 0.05,\, 0.10,\, 0.15,\, 0.20,\, 0.25,\, 0.30\}$ and $\ts$ is a directional vector considered as either $\ts=(1,1,\dots,1)^\top$ or $\ts=(1,0,\dots,0)^\top$. Notice that both distributions (a) and (b) use the identity scale matrix which allows us to compare easier the impact of shifts in the two different directions. 

The dimension was set as $d\in\{2,4,6\}$.
The sample size was chosen as $n\in\{150,300\}$ and the grid was constructed with $n_R=6$, so $\| \Fs(\tX_i)\|$ can take $n_R$ different values. %Note that a smaller sample size is not considered, because for dimension $d=6$ $n\leq 100$ is insufficient fo
The choice $n_R=15$ was considered as well, but the conclusion regarding the power comparison are analogous to those for $n_R=6$, so these results are not presented for the sake of brevity.

For the \emph{one-sample c-o test with random signs} from Section \ref{sec:test1}, the  
 grid $\mathcal{G}_n$ was constructed as one of the following:
\begin{description} 
\item[(R1)] regular random grid with $n_R$ spheres and $n_S=n/n_R$ directions which are generated randomly from the uniform distribution on $\mathcal{S}^{d-1}$,
\item[(R2)] random grid with $n_R$ spheres and  $n_S=n/n_R$ directions, which are drawn randomly from $\mathcal{S}^{d-1}$ for each sphere separately (so the resulting grid does not satisfy the factorization $\boldsymbol{g}=r_i \boldsymbol{u}_j$ mentioned in Section~\ref{sec:prel}), % (i.e.~generally $n$ different directions are generated in total), 
\item[(H)] grid with  $n_R$ spheres  and $n_S$ directions which are constructed from a Halton sequence of $n_S$ points in $\R^{d-1}$.
%\item[(S)] grid with  $n_R$ spheres  and $n_S$ directions which are constructed from a Sobol sequence of $n_S$ points in $\R^{d-1}$.
\end{description}
For a given $n_S$ and dimension $d$ the Halton sequence of points in $\mathbb{R}^{d-1}$ was generated using function \texttt{halton} from package \texttt{randtoolbox}, see \cite{randtool}.   Subsequently these points were transformed to $\mathcal{S}^{d-1}$ using polar coordinates and the corresponding suitable transformation described in \cite[Section 1.5.3]{fang}.  
The empirical c-o distribution function $\Fs$  was computed using so called Hungarian method with the help of  function \verb"solve_LSAP" from package \texttt{clue}, see \cite{clue}.

%Technical details for the construction of the latter two grids and the computation of $\Fs$ are provided  in Appendix~\ref{ap:A}. 
The  \emph{symmetrized one-sample} c-o test from Section~\ref{sec:test2} requires  a symmetric grid of $2n$ points satisfying \eqref{eq:Gsym}. It was constructed as
\begin{description}
\item[(R2*)] union of $\mathcal{G}$ and $-\mathcal{G}$, where $\mathcal{G}$ is obtained by (R2),
\item[(H*)] grid with  $n_R$ spheres, $n_S$ directions $\boldsymbol{s}$ which are constructed from a Halton sequence of $n_S$ points in $\R^{d-1}$ transformed to  $\mathcal{S}^{d-1}\cap \left(\R^{d-1}\times [0,\infty)\right)$, and corresponding  $n_S$ directions $-\boldsymbol{s}$.
%\item[(S*)]  analogous to (H*) but taking a Sobol sequence of $n_S$ points in $\R^{d-1}$.
\end{description}

All tests were computed for Wilcoxon scores with $J(u)=u$.  For some other possible choices for $J$ see \cite[Section~6]{marc2}, in particular their Section~6.2 for the discussion of the choice of the score function and efficiency of the test, and for a comparison in the context of two-sample tests see \cite{mordant}.
The power of the c-o rank tests is compared to the benchmark 
Hotelling's one-sample $T^2$ test and also to the one-sample spatial rank test from  \cite{OjaRan04},  see \cite{SpatialNP}.
%on hyperplane based signs and ranks. The latter concept of so called interdirections and lift-interdirections was proposed by \cite{OjaPain}, inspired by previous work of \cite{randles}. It proved to be a useful tool for elliptical distributions  in various situations, see also \cite{HaPain02a} and  \cite{HaPain02b}. 
%The computation of the original signs and ranks is very time demanding (or even impossible) for larger $d$, but recently \cite{HuSiKlic} and \cite{HuSi} recommended to use tests based on a randomized version of  interdirections and lift-interdirections,  and showed that the corresponding tests have comparable properties while they can be computed in a substantially shorter time. 
%Thus,  the one-sample test based on this randomized concept with $5\cdot n$ hyperplanes taken for computation of the ranks and signs and the Wilcoxon scores was computed in this Monte Carlo study.

To sum up, there are 7 tests  to be compared: 3 based on the  C-O test with random signs and with various grids (abbreviated as RAN-R1, RAN-R2, RAN-H), 2 based on the symmetrized sample and various grids (SYM-R2, SYM-H), the Hotelling's $T^2$ test (HOT) and the spatial rank test (SPAT).  The power of all tests is computed from $500$ replications.

\begin{figure}[htbp]
\centering
\includegraphics[width=\textwidth]{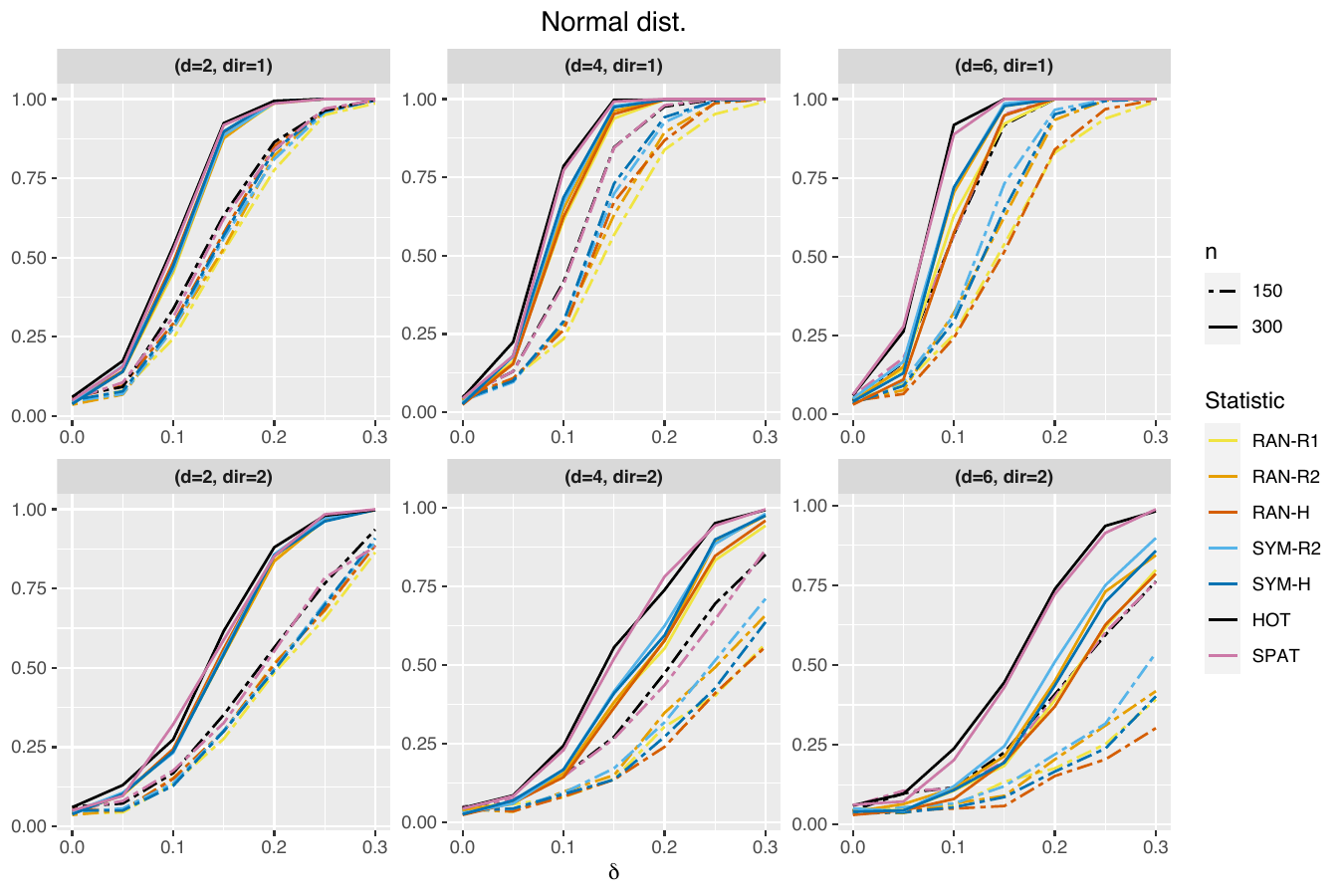}
\caption{Power of the considered tests for the standard normal distribution $\mathsf{N}_d(\boldsymbol{0},\boldsymbol{I})$ in $\R^d$ for $d\in\{2,4,6\}$ (in columns) shifted by $\delta \cdot \ts$ for $\ts=(1,\dots,1)^\top$ (dir$=1$, first row) or $\ts=(1,0,\dots,0)^\top$ (dir$=2$, second row) for  sample sizes $n\in\{150,300\}$ and various $\delta\geq 0$.}\label{fig:norm}
\end{figure}

\subsection{Results}

\begin{figure}[htbp]
\centering
\includegraphics[width=\textwidth]{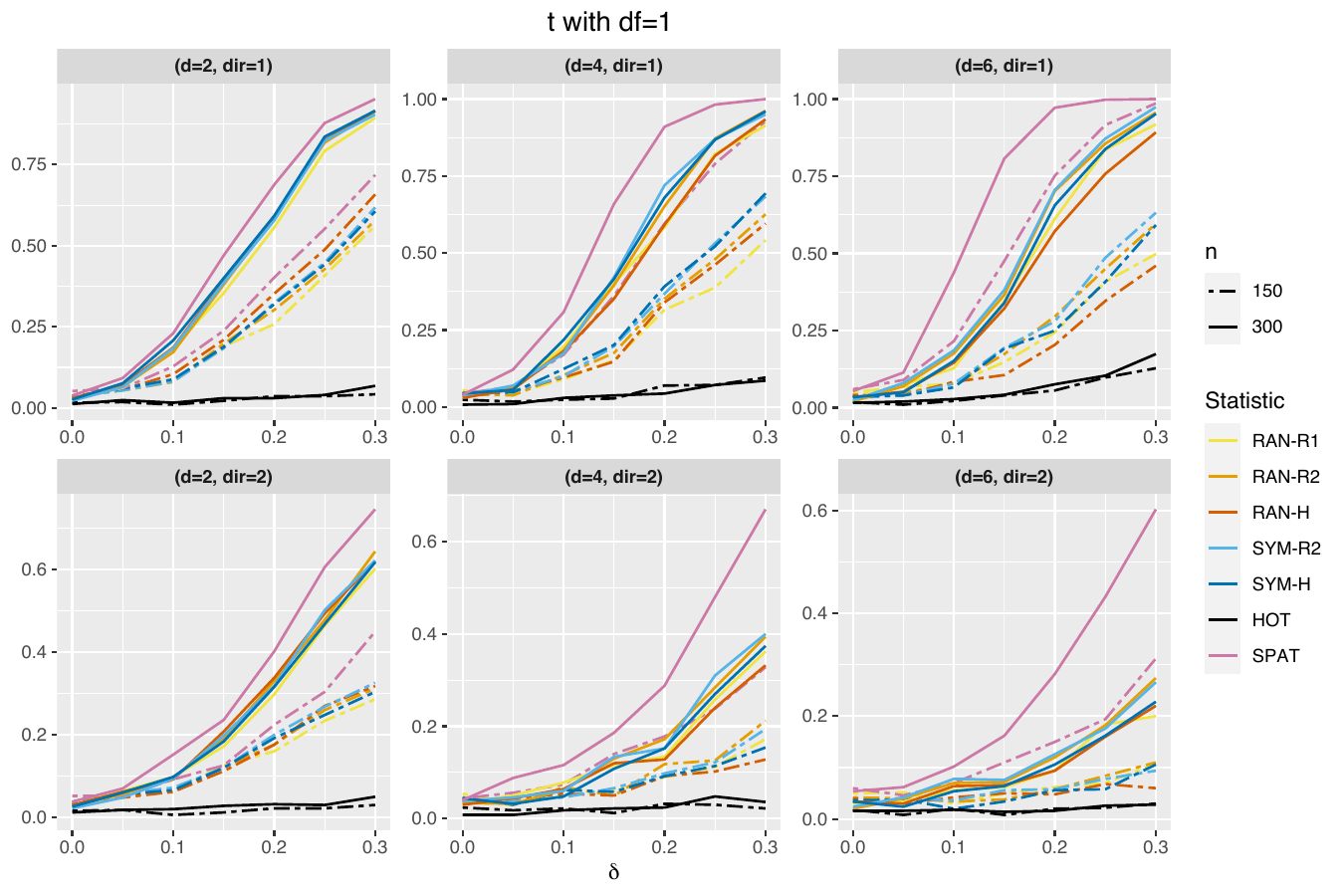}
\caption{Power of the considered tests for $t_1$ distribution in $\R^d$ for $d\in\{2,4,6\}$ (in columns) shifted by $\delta \cdot \ts$ for $\ts=(1,\dots,1)^\top$ (dir$=1$, first row) or $\ts=(1,0,\dots,0)^\top$ (dir$=2$, second row) for  sample sizes  $n\in\{150,300\}$ and various $\delta\geq 0$.}\label{fig:t.1}
\end{figure}

\begin{figure}[htbp]
\centering
\includegraphics[width=\textwidth]{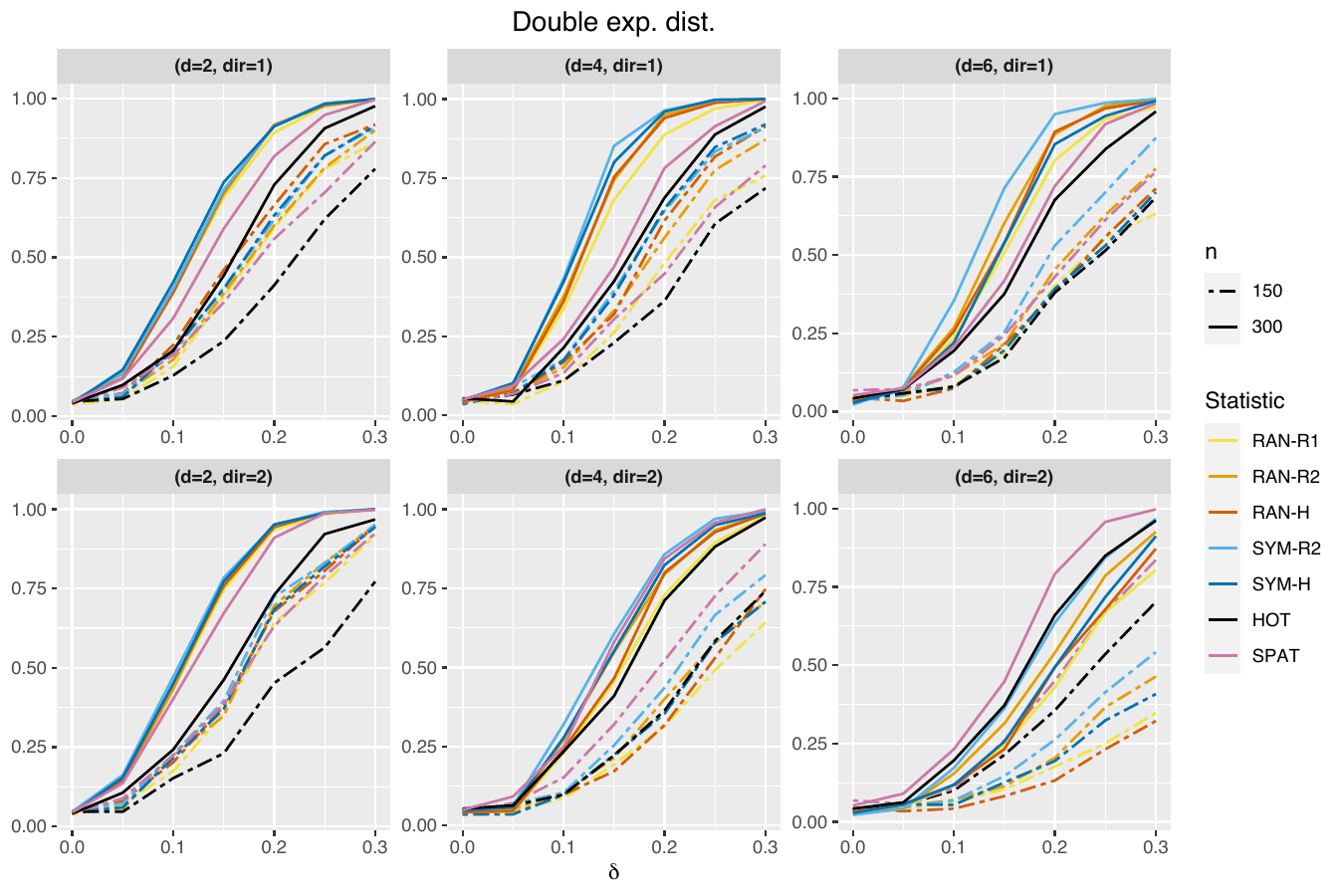}
\caption{Power for the ``double exponential'' distribution in $\R^d$ for $d\in\{2,4,6\}$ (in columns) shifted by $\delta \cdot \ts$ for $\ts=(1,\dots,1)^\top$ (dir$=1$, first row) or $\ts=(1,0,\dots,0)^\top$ (dir$=2$, second row) for  sample sizes  $n\in\{150,300\}$ and various $\delta\geq 0$.}\label{fig:exp}
\end{figure}

The obtained results are provided
 in Figures  \ref{fig:norm}--\ref{fig:exp}  in forms of power functions  depending on the size of the shift $\delta$ for the significance level $\alpha=0.05$. 
In each plot, the various considered tests are distinguished by colors and the sample sizes by line types (dotted-dashed for $n=150$, and solid for $n=300$).  The three considered dimensions $d=2,4,6$ are provided in the three plot columns, while the rows  correspond to the  two directional vector $\ts$,  with direction 1 and 2 corresponding to $\ts=(1,1,\dots,1)^\top$ and $\ts=(1,0,\dots,0)^\top$ respectively. In addition, the size of the tests is summarized numerically in Table~\ref{tab:size}. 

All considered c-o tests satisfy the prescribed significance level $\alpha=0.05$ with a tendency to be slightly undersized (conservative) in some settings. 
 Regarding the power, it is not surprising that the Hotelling's $T^2$ test leads to the largest power for the normal distribution. In this case, the power of c-o tests is generally comparable for $d=2$ and for $d=4$ with shift in direction $\ts=(1,1,\dots,1)^\top$. 
For $d=6$ larger differences are observed, but they decrease with an increasing sample size $n$. Among the c-o tests, the symmetrized version with the randomized grid (R2)  performs the best. 
The Hotelling's $T^2$ test completely fails for $t_1$ distribution, while the test based on spatial ranks performs the best. Differences in power between this test and the c-o tests seem to increase with a growing dimension $d$. 
Finally, the double exponential distribution illustrates the benefits of c-o ranks based tests, because in this case %the symmetrized 
all the proposed c-o tests  completely outperform the Hotelling's $T^2$ test as well as the spatial rank tests if the shift is in direction  $\ts=(1,1,\dots,1)^\top$. For shift in direction  $\ts=(1,0,\dots,0)^\top$, the symmetrized test performs the best for dimension $d\leq 4$, but it is outperformed by the spatial rank test for $d=6$.

Comparing the two proposed c-o approaches (random signs and symmetrized), it is not possible to say that the symmetrized test leads to a larger power compared to the randomized test in all cases, but it could be slightly preferred. Its disadvantage, compared to the randomized test, is that it uses a transport of $2n$ points, so it is more computationally demanding for larger sample sizes $n$. 
 Regarding the choice of the grid, for both randomized and symmetrized test, the choice of the random grid (R2) seems to lead to overall the best results compared to a random grid (R1) and to the fixed regular grid (H).

\begin{table}[ht]
\centering
\begin{tabular}{rr|rrrrrrr}
  \hline
\hline
$d$ & $n$ & RAN-R1 & RAN-R2 & RAN-H & SYM-R2 & SYM-H & HOT & SPAT \\ 
\hline
& &\multicolumn{7}{c}{Normal dist.}\\
   \hline
 2 & 150 & 0.040 & 0.036 & 0.056 & 0.046 & 0.050 & 0.062 & 0.052 \\ 
   & 300 & 0.044 & 0.040 & 0.044 & 0.044 & 0.038 & 0.060 & 0.050 \\ 
  4 & 150 & 0.036 & 0.042 & 0.048 & 0.036 & 0.040 & 0.050 & 0.048 \\ 
   & 300 & 0.040 & 0.040 & 0.024 & 0.030 & 0.026 & 0.046 & 0.044 \\ 
  6 & 150 & 0.044 & 0.040 & 0.042 & 0.034 & 0.034 & 0.044 & 0.058 \\ 
   & 300 & 0.050 & 0.040 & 0.030 & 0.050 & 0.040 & 0.060 & 0.062 \\ 
   \hline
& &\multicolumn{7}{c}{$t_1$ dist.}\\
   \hline
2 & 150 & 0.030 & 0.038 & 0.036 & 0.034 & 0.030 & 0.016 & 0.052 \\ 
   & 300 & 0.028 & 0.028 & 0.030 & 0.022 & 0.026 & 0.012 & 0.038 \\ 
  4 & 150 & 0.054 & 0.042 & 0.038 & 0.036 & 0.046 & 0.024 & 0.044 \\ 
   & 300 & 0.034 & 0.034 & 0.030 & 0.038 & 0.044 & 0.008 & 0.040 \\ 
  6 & 150 & 0.050 & 0.042 & 0.040 & 0.032 & 0.036 & 0.018 & 0.060 \\ 
   & 300 & 0.030 & 0.020 & 0.034 & 0.024 & 0.034 & 0.016 & 0.054 \\ 
 \hline
& &\multicolumn{7}{c}{Double exp. dist.}\\
   \hline
2 & 150 & 0.036 & 0.044 & 0.048 & 0.048 & 0.044 & 0.046 & 0.042 \\ 
   & 300 & 0.042 & 0.046 & 0.046 & 0.046 & 0.040 & 0.040 & 0.044 \\ 
  4 & 150 & 0.046 & 0.040 & 0.042 & 0.034 & 0.036 & 0.046 & 0.040 \\ 
   & 300 & 0.044 & 0.050 & 0.042 & 0.050 & 0.048 & 0.054 & 0.050 \\ 
  6 & 150 & 0.040 & 0.038 & 0.044 & 0.034 & 0.034 & 0.040 & 0.068 \\ 
   & 300 & 0.032 & 0.038 & 0.032 & 0.022 & 0.028 & 0.042 & 0.052 \\ 
   \hline
\hline
\end{tabular}
   \caption{Empirical size of the considered tests for  the three $d$-dimensional distributions:
   the randomized one-sample test with various grids (RAN-R1, RAN-R2, and RAN-H), the symmetrized test with two different grids (SYM-R2 and SYM-H), Hotelling's $T^2$ test (HOT), and the spatial rank test (SPAT).
   }\label{tab:size}
\end{table}

\subsection{Sensitivity to the assumption of central symmetry}

The proposed c-o tests rely on the assumption of central symmetry of the true distribution. One may ask, how the tests perform if this assumption is violated.  Of course, if the distribution of $\tX$ is not centrally symmetric, then there is no true point $\tmu_S$ to be tested, and, in fact, it is not clear what is the theoretical counterpart of the location estimator $(F_{\pm}^{(n)})^{-1}(\boldsymbol{0})$. However, it still can be of interest to explore the behavior of the proposed test for small deviations from the symmetry.

%In order to explore their sensitivity  with respect to this assumption, 
%the empirical size and power  are computed for samples simulated from a bivariate skew normal distribution.  We make use of 
For this purpose, we use a bivariate skew normal distribution $\mathsf{SN}_2(\boldsymbol{0},\boldsymbol{I},\boldsymbol{\alpha})$ for $\boldsymbol{\alpha}=(\alpha,\alpha)^\top$, $\alpha\in\R$, where the parametrization is taken from \cite[Chapter~5]{SN_book}, i.e. the density of this distribution is
\[
f(\boldsymbol{x}) = 2 \varphi(\boldsymbol{x})\Phi(\boldsymbol{\alpha}^\top\boldsymbol{x}),
\]
where $\varphi$ is the density of a bivariate normal distribution $\mathsf{N}_2(\boldsymbol{0},\boldsymbol{I})$ and $\Phi$ is the cumulative distribution function of univariate $\mathsf{N}(0,1)$. The parameter $\boldsymbol{\alpha}$ is referred to as a slant parameter, and it drives the skewness of the distribution. The mean of this distribution is $\boldsymbol{\mu}=(\mu,\mu)^\top$ for $\mu=\sqrt{2/\pi}\alpha/\sqrt{1+2\alpha^2}$.
%so under the null hypothesis we simulate data as $\tX=\boldsymbol{Y}-\boldsymbol{\mu}$, where $\boldsymbol{Y}\sim \mathsf{SN}(\boldsymbol{0},\boldsymbol{I},\boldsymbol{\alpha})$. Under an alternative, we take 

\begin{figure}[htbp]
\centering
\includegraphics[width=\textwidth]{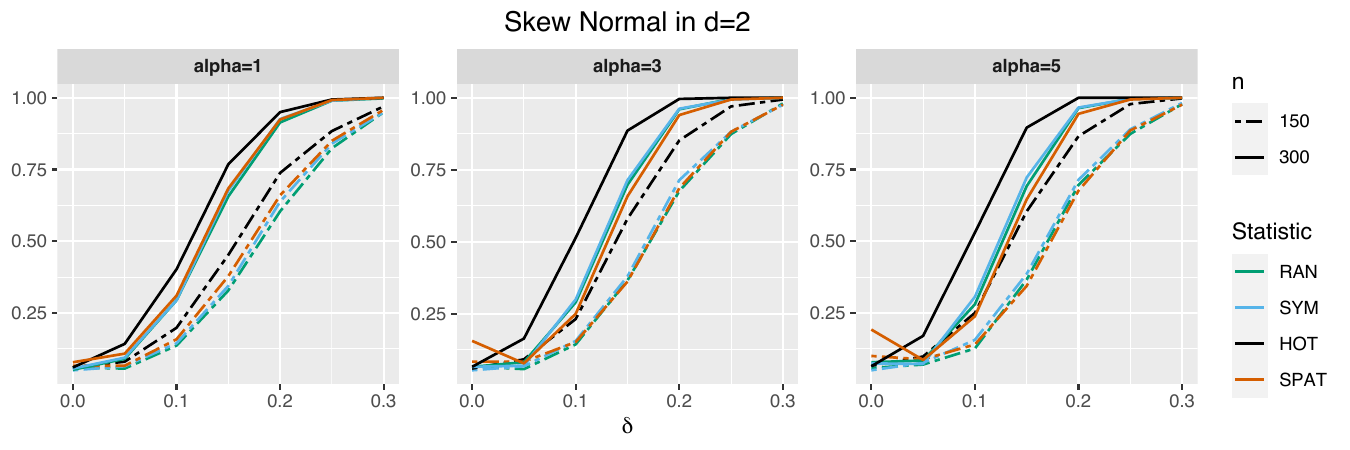}
\caption{Proportion of rejections of the null hypothesis of zero location for the c-o test with random signs (RAN), symmetrized sample (SYM), Hotelling's $T^2$ test (HOT), and spatial rank test (SPAT) for skew data in dimension $d=2$. 
%Under the alternative, the distribution is shifted by $\delta\ts$ for 
% $\ts=(1,1)^\top$.
The columns correspond to distributional parameter $\alpha$: the larger $\alpha$, the more skewed distribution, see details in the main text.}\label{fig:SN}
\end{figure}

 We simulate data  as $\tX=\boldsymbol{Y}-\boldsymbol{\mu}+\delta\cdot \boldsymbol{s}$ for $\boldsymbol{Y}\sim \mathsf{SN}(\boldsymbol{0},\boldsymbol{I},\boldsymbol{\alpha})$ and
 a directional vector $\boldsymbol{s}=c(1,0)$ and a shift size $\delta\geq 0$. 
 For each sample, we computed the following tests of zero location:  test based on random signs (RAN, grid (R2)),  test based on a symmetrized sample (SYM, grid (R2)), the Hotelling's $T^2$ test (HOT) and the spatial rank test (SPAT). Here, ``location'' is the mean for the Hotelling's $T^2$ test, while it is more complicated to define the theoretical location for three remaining tests for skewed data.
 Subsequently, proportions of rejections for significance level $0.05$ were calculated and they are plotted in Figure~\ref{fig:SN} for  $n\in\{150,300\}$ and for $\alpha\in\{1,3,5\}$.

 For $\delta=0$ and $\alpha=1$, all the four procedures 
reject the null hypothesis of zero location in approximately 5 \% of the cases, while for $\delta=0$ and $\alpha>1$ this is the case only for the first three tests. This is not surprising for the Hotelling's $T^2$ test, because for $\delta=0$ we have $\E\tX=\boldsymbol{0}$ regardless the value of the slant parameter $\alpha$, so the null hypothesis holds and the test keeps the nominal level asymptotically. In addition,  $\E\tX\ne \boldsymbol{0}$ for $\delta>0$, the test is consistent, and
it leads to the largest percentage of rejections for any fixed $\delta>0$.
It is interesting to see that the two proposed c-o tests (RAN and SYM) seem to be quire robust with respect to small deviations from the symmetry as they also reject the null in approximately 5 \% cases for $\delta=0$ and the proportion of rejections growths with $\delta$ in an expected way.

 %The proportion of rejections 
 %$n\in\{150,300\}$ and for $\alpha\in\{1,3,5\}$ are provided in Table~\ref{tab:SN}. For simplicity, only one version of the test with random signs (RAN, grid (R2)) and one version of the test based on a symmetrized sample (SYM, grid (R2)) are compared to the Hotelling's $T^2$ test (HOT) and spatial rank test (SPAT). 
 %Clearly, the skewness of the distribution most strongly affects the size of the spatial rank test. For $\alpha=5$ the nominal size 5 \% is exceeded more than three times. The empirical sizes of c-o tests and the Hotelling's $T^2$ test are comparable, all with a tendency to be slightly larger than the nominal level 5 \%, in particular for $\alpha=5$ and $n=300$, but the size is never larger than 8 \%. Hence, we may conclude that a small deviation from the central symmetry does not affect the size of the proposed c-o one-sample tests dramatically.   

%The power of the four considered tests for shift in direction  $\ts=(1,0)^\top$ is compared  in Figure~\ref{fig:SN}. Note that the results for direction $\ts=(1,1)^\top$ are analogous, so they are not presented. 
%The Hotelling's $T^2$ test leads overall to the largest power, but the c-o tests behave reasonably well, even for $\alpha=5$. 
%The spatial rank test, based on the ellipticity assumption, exceeds the test size, so discussion regarding its power is irrelevant. 

\subsection{Comparison with marginal one-sample test}

\begin{figure}[htbp]
\centering
\includegraphics[width=\textwidth]{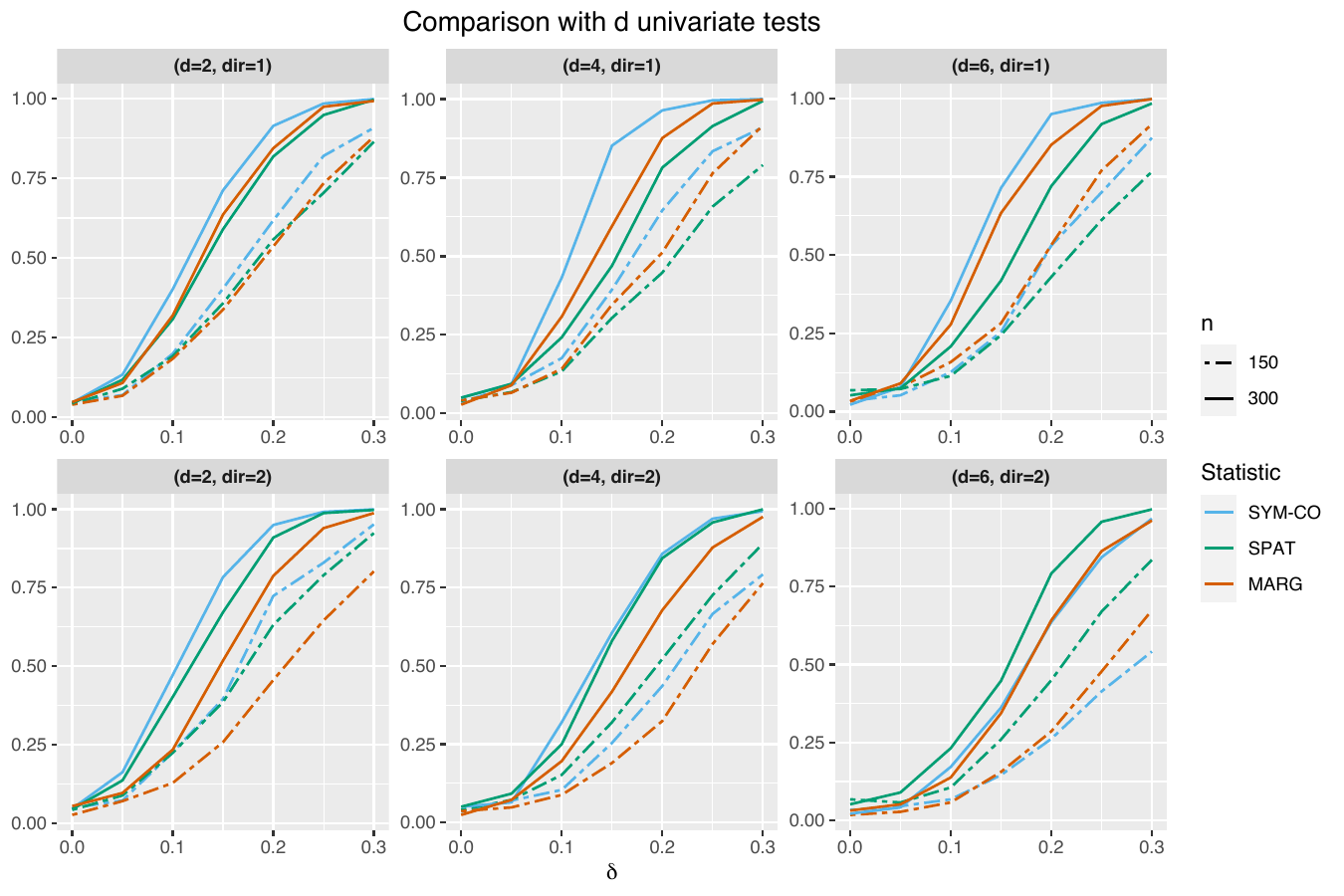}
\caption{Power curves of two multivariate one-sample tests (symmetrized c-o test and spatial rank test) compared to the power of $d$ univariate marginal Wilcoxon tests (MARG) combined via the Bonferroni correction for  ``double exponential'' distribution in $\R^d$ for $d\in\{2,4,6\}$ (in columns). A centered sample is shifted by $\delta \cdot \ts$ for $\ts=(1,\dots,1)^\top$ (dir$=1$, first row) or $\ts=(1,0,\dots,0)^\top$ (dir$=2$, second row) for  sample sizes  $n\in\{150,300\}$.
}\label{fig:marg}
\end{figure}

One of the reviewers raised an interesting question regarding the benefits of the multivariate one-sample test with respect to $d$ univariate tests with a correction for multiple testing. The corresponding power comparison is provided in Figure~\ref{fig:marg}, where the power of the symmetrized c-o test with a random grid (R2) and the power of the spatial rank test is compared to the power of $d$ univariate Wilcoxon tests (conducted for each coordinate) combined via the Bonferroni correction. Results for the double exponential distribution are presented, but the overall conclusion is the same also for the normal and $t_1$ distribution.  The results reveal that the comparison depends on the dimension $d$ and on the shift direction, while differences between sample sizes $n$ are not so visible. Indeed, there are situations (shift in a suitable direction) where it is more beneficial to simply conduct $d$ univariate marginal tests instead of one multivariate test (spatial or c-o test), but the opposite holds for some other settings.  In practice, one typically does not have a prior information about the direction of the possible shift, so a multivariate test is a more general choice. 

 % \section{Complements to the simulation study}\label{ap:A}

%\begin{table}[ht]
%\centering
%\begin{tabular}{rr|rrrr}
%  \hline  \hline
%$n$ & $\alpha$ & RAN & SYM & HOT & SPAT \\ 
%  \hline
%150 & 1 & 0.060 & 0.050 & 0.060 & 0.066 \\ 
%   & 3 & 0.068 & 0.054 & 0.060 & 0.084 \\ 
%   & 5 & 0.058 & 0.050 & 0.066 & 0.100 \\ 
%  300 & 1 & 0.052 & 0.056 & 0.060 & 0.078 \\ 
%   & 3 & 0.066 & 0.066 & 0.066 & 0.156 \\ 
%   & 5 & 0.078 & 0.074 & 0.064 & 0.192 \\ 
%   \hline  \hline
%\end{tabular}
%\caption{Empirical size of a test based on random signs (RAN), test %based on symmetrized sample (SYM), Hotelling's $T^2$ test (HOT), and a %spatial rank test (SPAT) for data simulated from a bivariate skew normal %distribution with zero mean. The parameter $\alpha$ drives the skewness: %the larger $\alpha$, the more skewed distribution, see details in the %main text.}\label{tab:SN}
%\end{table}

 \bibliographystyle{apalike}
 \bibliography{CO-lit}

\end{document}